\documentstyle[12pt]{article}
\newenvironment{beweis}{{\it Proof.}\ }{\ $\ \ \ \Diamond$ \\\ }

 \newcounter{nsatz}[section]
 \newcounter{nlemma}[section]
 \newcounter{ndef}[section]
 \newcounter{nhyp}[section]
 \newcounter{nconjecture}[section]
 \newcounter{ncor}[section]
 \newcounter{nrem}[section]
 \newcounter{nexample}[section]
 \newcounter{nprop}[section]

 \newenvironment{nsatz}{\refstepcounter{nsatz}{\bf \arabic{section}.\arabic{nsatz}}\
                {\sc\bf Theorem.\ }\it}{\\\\ \rm}

 \newenvironment{nlemma}{\setcounter{nlemma}{\value{nsatz}}
                \refstepcounter{nlemma}
                \setcounter{nsatz}{\value{nlemma}}
                {\bf \arabic{section}.\arabic{nsatz}}\
                {\sc\bf Lemma.\ }\it}{\\\\ \rm}

 \newenvironment{ncor}{\setcounter{ncor}{\value{nsatz}}
                \refstepcounter{ncor}
                \setcounter{nsatz}{\value{ncor}}
                {\bf \arabic{section}.\arabic{nsatz}}\
                {\sc\bf Corollary.\ }\it}{\\\\ \rm}

\begin{document}
\newcommand{\n}{{\mbox{\rm I$\!$N}}}
\newcommand{\z}{{\mbox{{\sf Z\hspace{-0.4em}Z}}}}
\newcommand{\R}{{\mbox{\rm I$\!$R}}}
\newcommand{\Q}{{\mbox{\rm I$\!\!\!$Q}}}
\newcommand{\C}{{\mbox{\rm I$\!\!\!$C}}}
\newcommand{\p}{{\mbox{\rm I$\!\!\!$P}}}
\newcommand{\ug}{\ \raisebox{-.3em}{$\stackrel{\scriptstyle \leq}
{\scriptstyle \sim}$} \ }

\thispagestyle{empty}
\setlength{\parindent}{0pt}
\setlength{\parskip}{5pt plus 2pt minus 1pt}

%\frenchspacing
%\sloppy
%\thispagestyle{empty}{\vspace*{-2cm} \hspace*{-3.5cm}
\thispagestyle{empty}
\newcommand{\Irr}{{\mbox{\rm Irr}}}
\newcommand{\sIrr}{{\mbox{\scriptsize\rm Irr}}}
\newcommand{\Char}{{\mbox{\rm Char}}}
\mbox{\vspace{4cm}}
\vspace{4cm}
\begin{center}
{\bf \Large\bf  A lower bound for the number of conjugacy classes of finite groups\\}
\vspace{3cm}
by\\
\vspace{11pt}
Thomas Michael Keller\\
Department of Mathematics\\
Texas State University\\
601 University Drive\\
San Marcos, TX 78666\\
USA\\
e--mail: keller@txstate.edu\\
\vspace{1cm}
\end{center}
%Short form of title: 'Orbit sizes of Permutation groups'.
\thispagestyle{empty}
\newpage

\begin{center}
\parbox{12.5cm}{{\small
{\sc Abstract.}
In 2000, L. H\'{e}thelyi and B. K\"{u}lshammer proved that if $p$ is a prime number dividing the order of a finite
solvable group $G$, then $G$ has at least $2\sqrt{p-1}$ conjugacy classes. In this paper we show that if $p$ is large,
the result remains true for arbitrary finite groups.\\
}}
\end{center}
\normalsize

\section{Introduction}\label{section0}
Let $G$ be a finite group and write $k(G)$ for the number of its conjugacy classes. Finding general bounds on
$k(G)$ is a fundamental problem in finite group theory, and in this paper we are concerned with lower bounds. There
is already a large body of work on this topic, and one of the strongest general results is due to Pyber \cite{pyber}
and states that there is an $\epsilon>0$ such that
\[k(G)\geq\epsilon\frac{\log|G|}{\log\log|G|}\]
Many other lower bounds on $k(G)$ in more specialized situations can be found in
\cite{bertram}.\\
Recently, L. H\'{e}thelyi and B. K\"{u}lshammer added a new twist to the subject. In \cite{hethelyi} they proved
that if $G$ is solvable and $p$ is a prime dividing $|G|$, then
\[k(G)\ \geq\ 2\sqrt{p-1}.\]
In this paper we show that this result remains true for arbitrary groups if $p$ is large. More precisely, we
prove\\

{\bf Theorem A. }{\it
There is a constant $C$ such that if $G$ is a finite group whose order is divisible by a prime $p>C$,
then $k(G)\geq 2\sqrt{p-1}$.
}\\

This is \ref{cor8} below.\\
We remark that we do not use \cite{hethelyi} in our proof of Theorem A. We do, however, heavily rely on the results
recently obtained by G. Malle \cite{malle}, and as such the proof relies on the Classification of Finite Simple Groups
(and also on GAP).\\

Thanks to a reduction result in \cite{malle} the proof of Theorem A reduces to proving the following
result.\\

{\bf Theorem B. }{\it
With $C$ as in Theorem A, we have the following.\\
Let $G$ be a finite group and let $V$ be a finite, faithful, irreducible $GF(p)G$--module with $p$ a prime not
dividing $|G|$. If $p>C$, then
\[k(GV)\ \geq\ 2\sqrt{p-1}.\]
}\\

This is \ref{theo7} below.\\
Note that the situation studied here is that of the well--known $k(GV)$--problem (see e. g. \cite{keller}) whose
goal it is to find an upper bound for $k(GV)$ under the hypothesis of Theorem B (namely, $k(GV)\leq |V|$). So
here we approach $k(GV)$ from the other side, seeking a general
lower bound.\\

Perhaps unexpectedly (and unlike the proofs of upper bounds for $k(GV)$), our proof of Theorem B is not
inductive.\\

We also remark that while it was shown in \cite{hethelyi} that $k(G)=2\sqrt{p-1}$ is possible for suitable $G$,
it seems that the bound $k(G)\geq 2\sqrt{p-1}$ is fairly weak in general, and the proof of Theorem B gives some
further evidence of this.\\

Notation: All logarithms are to base 2 throughout the paper. $F^*(G)$ denotes the generalized Fitting subgroup of the
group $G$. If $G$ acts on a set $\Omega$, then $n(G,\Omega)$ denotes the number of orbits of $G$ on $\Omega$.
$S_n$ denotes the symmetric group on $n$ letters, and $A_n$ is the alternating group. $A\ug B$ means that the group $A$
is isomorphic to a subgroup of the group $B$. If $V$ is a $G$--module and $H\leq G$, then $V_H$ denotes $V$ considered
as an $H$--module. By $V(d,q)$ ($q$ a prime power) we denote $GF(q)^d$, the $d$--dimensional vector space over
$GF(q)$.
If $G\ug S_n$ and $H$ is a group, then by $H\wr G$ we denote the corresponding wreath product of $H$ and $G$.
All other notation is standard.\\

\section{Results and proofs}\label{section1}

Before we start we recall the well--known result that if $G$ is a group and $H\leq G$, then
\[k(G)\geq\frac{k(H)}{|G:H|}\]
We will use this fact freely throughout the paper. \\
We first use a result of Malle \cite{malle} in the following
version.\\

\begin{nsatz}\label{theo1}
Let $p$ be a prime and $q$ be a power of $p$. Let $G$ be a finite group with quasisimple generalized Fitting subgroup
$F^*(G)$, and let $V\cong V(d,q)$ be an absolutely irreducible, faithful $GF(q)G$--module with $(q,|G|)=1$. Then
one of the following holds:\\

(i) $d\leq 12$\\

(ii) $n(G,\p(V))\geq p$, where $\p(V)$ denotes the projective space of $V$.
\end{nsatz}
\begin{beweis}
Let $K\subseteq GF(q)$ be the smallest field over which the representation of $G$ on $V$ can be realized, i. e.
choose $K$ minimal such that there is a $KG$--module $W$ with $V=W\oplus_K GF(q)$.
Then $W$ is an absolutely irreducible $KG$--module that cannot be realized over a proper subfield of $K$.
Now suppose that $d>12$.
By \cite[Satz 5.1]{malle} it follows that $n(G,\p(W))\geq p$. Thus clearly $n(G,\p(V))\geq p$, as wanted.
\end{beweis}

\begin{ncor}\label{cor2}
There is a universal constant $C_0$, such that the following holds. Let $G$ be a finite group and $V$ be a finite faithful
$G$--module of characteristic $p$ such that $(p,|G|)=1$. Suppose that $N$ is a normal subgroup of $G$ which is a central
product of a quasisimple group $G$, and the cyclic group $Z=Z(F^*(G))$; write $N=G_1\circ Z$ for this central product.
Furthermore assume that $V_N$ is homogeneous. If $|G_1|>C_0$, then
\[n(G,V)\geq\frac{p}{|G:N|}.\]
\end{ncor}
\begin{beweis}
Clearly it suffices to show that $n(N,V)\geq p$. Now let $V_0\leq V$ be an irreducible $N$--submodule of $V$. As
$V_N$ is homogeneous, $V_0$ is also faithful as $N$--module. Then it suffices to show that
$n(N,V_0)\geq p$.\\
Now by \cite[Lemma 10]{robinson-thompson} and \cite[Theorem 26.6]{aschbacher} it is clear that if we put
$K={\rm End}_{GF(p)G}(V_0)$, then $K=GF(q)$ for some power $q$ of $p$, and if $W$ is an irreducible $N$--submodule
of the $N$--module $V_0\oplus_{GF(p)}K$, then $W$ is an absolutely irreducible $KN$--module (which can not be
realized over a proper subfield of $K$), and the permutation actions of $G$ on $V_0$ and $W$ are equivalent. Hence
$n(N,V_0)=n(N,W)$.\\
Now as $W$ is absolutely irreducible, clearly $W$ is irreducible as $G_1$--module (as $Z$
acts on $W$ by scalar multiplication) and $p$ does not divide $|G|$, by Jordan's Theorem (see
\cite[Theorem 15.7]{dixon}) we see that we can choose $C_0$ large enough so that $|G_1|>C_0$ implies
$\dim W>12$. thus by \ref{theo1} we see that $n(G_1,\p(W))\geq p$, and since $Z$ acts trivially on $\p(W)$, it
follows that $n(N,\p(W))\geq p$. Hence $n(N,W)\geq p$, and so we are done.
\end{beweis}

Next we need to recall the following result by Gambini--Weigel and Weigel \cite{gambini}. We state it as in
\cite[Theorem 2.1]{gluck/magaard}, but take the opportunity to add a recent correction to
it.\\

\begin{nsatz}\label{theo4}
Let $G$ be a finite group and $W$ be a faithful primitive finite $G$--module with $(|G|,|W|)=1$. Then
\[|G|\leq |W|\log |W|,\]
except in the following cases:\\

(i) $|W|=7^4$ and $G$ is ${\rm Sp}(4,3)$ or $Z_3\times {\rm
Sp}(4,3)$\\

(ii) $|W|=3^4$, $|G|=4\cdot 5\cdot 2^5$, $G$ has exactly two orbits on $W$, and $G$ is the group $G_{3^4,1}^0$,
in \cite[Hauptsatz]{huppert1957}. \\
\end{nsatz}

The case (ii) was omitted in \cite{gambini}, but when informed about this omission, T. S. Weigel confirmed that this
would be the only other exception \cite{weigel}.\\

We only need the following consequence:\\

\begin{ncor}\label{cor5}
Let $G$ be a finite group and $W$ be a faithful primitive finite $G$--module with $(|G|,|W|)=1$. Then
\[|G|\leq 6|W|\log |W|\]
\end{ncor}

Next we need to study permutation groups with few conjugacy classes.\\

\begin{nlemma}\label{lem6}
There is a universal constant $A$ such that the following holds:
Let $G$ be a transitive permutation group on the set $\Omega$. Let $p$ be a prime and suppose that
$k(G)\leq 2\sqrt{p-1}$, and put $|\Omega|=n$.
If $p>A$, then
\[|G|\leq (\log p)^{8n}.\]
\end{nlemma}
\begin{beweis}
For $m\in\n$, let $p(m)$ be the number of partitions of $m$. It is well--known that $k(S_m)=p(m)$ and that
\[p(m)\ \sim\ \frac{e^{\pi\sqrt{\frac{2m}{3}}}}{4m\sqrt{3}}\qquad;\]
in particular, there exists a constant $C_1>0$ such that if $H=A_m$ or $S_m$,
then
\[(1)\qquad k(H)\geq C_1 2^{m^{\frac{1}{3}}}.\]
Let $G$ be a counterexample to the lemma with $n$ minimal.\\
Define $f:\ \n\longrightarrow \R$ by $f(m)=(\log p)^{8m}$ for $m\in\n$.
Now let $B\subseteq\Omega$ be a block (i. e., for each $g\in G$, $B^g\cap B=\emptyset$ or $B^g\cap B=B$)
which is minimal subject to $|B|>1$. (Thus $G$ is primitive if and only if $B=\Omega$.) It is then well--known
that $\Omega$ can be partitioned into subsets $B_i$, $i=1,\ldots,k$ for some $k\in\n$, such that $B=B_1$
and $G$ permutes the $B_i$ transitively. Let $K$ be the kernel of that action, so $G/K\ug S_k$, and let
$G_0$ be the setwise stabilizer in $G$ of $B$ and let $N_0\unlhd G_0$ be its pointwise stabilizer, and let
$\overline{G_0}=G_0/N_0$. Then
\[(2)\qquad G\ug\overline{G_0}\wr G/K,\]
as is well-known, and $\overline{G_0}$ is a primitive permutation group on
$B$.\\
Now if $K=1$, then $G\ug S_k$ with $k<n$, and hence by our choice of $G$ we get
\[|G|\leq f(k)\leq f(n),\]
and we are done.
So from now on we may assume that $K>1$. Hence clearly $k(G/K)\leq k(G)\leq 2\sqrt{p-1}$, and so by our minimal
choice of $K$ we have
\[(3)\qquad |G/K|\leq f(k)\]
Put $l=\frac{n}{k}$, so $|B|=l$. If $\overline{G_0}$ is not isomorphic to $A_l$ or $S_l$, then by
\cite[Corollary 1.2]{maroti} we know that $|\overline{G_0}|<3^l$ and thus by (2) and (3)
\[|G|\leq|G/K|\cdot |\overline{G_0}|^k\leq f(k)\cdot
3^{lk}=f(k)\cdot 3^n.\]
Next observe that as $|B|>1$, we have $k\leq\frac{n}{2}$, and so if $A$ is sufficiently large, we further obtain
$|G|\leq f(\frac{n}{2})\cdot 3^n\leq f(n)$, as wanted.

So now assume that $|\overline{G_0}|\cong A_l$ or $|\overline{G_0}|\cong S_l$, so $K\ug\left( S_l\right)^k$,
the direct product of $k$ copies of $S_l$.\\
Then it is clear that $G$ has a normal subgroup $N$ such that $N$ is the direct product of $t$ copies of $A_l$,
where $t$ is a suitable integer with $1\leq t\leq k$, and $|K/N|\leq 2^k$. Thus with (3) we get
\[(4)\qquad |G|=|G/K|\ |K/N|\ |N|\leq f(k) 2^k|A_l|^t\leq f(k)\cdot
2^k\cdot\left(l^l\right)^t\]
Next observe that by (3) and (1) we have
\[2\sqrt{p-1}\geq
k(G)\geq\frac{k(N)}{|G:N|}\geq\frac{k(A_l)^t}{|G:K||K:N|}\geq\frac{C_1^t 2^{tl^\frac{1}{3}}}{f(k) 2^k}\]
and hence $l\leq\left(\frac{1}{t}(k+1+\frac{1}{2}\log p+\log f(k))-\log
C_1\right)^3\leq\left(11\frac{k}{t}\log p-\log C_1\right)^3$.\\
Thus if $A$ is large enough, we see that this implies that there is a constant $C_2>1$ such that
\[(5)\qquad l\leq C_2\left(\frac{k}{t}\log p\right)^3\]
Using this and the fact that $l=\frac{n}{k}$, with (4) we find that
\begin{eqnarray*}
|G|&\leq &f(k) 2^k\left(C_2\left(\frac{k}{t}\log p\right)^3
         \right)^\frac{nt}{k}\\
   &\leq&f\left(\frac{n}{2}\right)\cdot
   2^\frac{n}{2}\left(\left(\frac{k}{t}\right)^\frac{t}{k}\right)^{3n}\left(C_2(\log
   p)^3\right)^\frac{nt}{k}
\end{eqnarray*}
Now clearly $\frac{t}{k}\leq 1$. Also as the function
$g(x)=x^\frac{1}{x}$ ($x>0$) is bounded above by 2, we conclude that
\begin{eqnarray*}
|G|&\leq &f\left(\frac{n}{2}\right)\cdot 2^{\frac{7}{2}n}\cdot\left(C_2\log
          p\right)^{3n}\\
   &\leq&\left(\log p\right)^{4n}\left(2^\frac{7}{6}C_2\log
   p\right)^{3n}\\
   &=&\left(2^\frac{7}{6}\cdot C_2\right)^{3n}\cdot (\log p)^{7n}
\end{eqnarray*}
Hence we have $|G|\leq f(n)$ if
\[(2^\frac{7}{6}C_2)^3\leq\log p,\]
but this is certainly true if $A$ is chosen sufficiently large. So the lemma is proved.
\end{beweis}

Now we can prove our main result.\\

\begin{nsatz}\label{theo7}
There is a universal constant $C$ such that the following holds.\\
Let $G$ be a finite group, $p$ be a prime not dividing $|G|$, and
let $V$ be a finite faithful, irreducible $GF(p)G$--module. If
$p>C$, then
\[k(GV)\geq 2\sqrt{p-1}.\]
\end{nsatz}
\begin{beweis}
Working towards a contradiction, assume that
\[k(GV)< 2\sqrt{p-1}.\]
Then clearly
\[2\sqrt{p-1}>k(GV)>n(G,V)\geq\frac{|V|}{|G|}\]
and thus
\[(6)\qquad |V|\leq 2\sqrt{p-1}\ |G|.\]
Put $m=\dim V$, so $|V|=p^m$.\\
If $m=1$, then $k(GV)=\frac{p-1}{|G|}+|G|$, and this easily implies the
assertion.\\

So now suppose $m\geq 2$.\\
It is well--known that one can find a subspace $V_1\leq V$ such that
if $H=N_G(V_1)$, then $H_1:=H/C_G(V_1)$ acts primitively and
faithfully on $V_1$, and $V$ is induced from the $H$--module $V_1$,
and $G$ transitively permutes the elements in $\Omega=\{V_1^g\ |\
g\in G\}$, and $V=V_1\oplus\ldots\oplus V_n$ for some $n\in\n$ and
suitable $V_i\in\Omega$ ($i=2,\ldots,n$), and
$\Omega=\{V_1,\ldots,V_n\}$. (Note that possibly $V_1=V$, and in this case $n=1$
and $G\cong H_1$ acts primitively on $V$.)\\
Let $K$ be the kernel of
the permutation action of $G$ on $\Omega$. Then $G/K\ug S_n$, and
clearly $k(G/K)\leq 2\sqrt{p-1}$. Then by \ref{lem6} we know that, if $C$ is large enough,
\[(7)\qquad |G/K|\leq\left(\log p\right)^{8n}\]
Moreover,
\[(8)\qquad G\ug H_1\wr (G/K)\]
and thus $K\ug H_1^n$ the direct product of $n$ copies of
$H_1$.\\
Now by (6), (7) and (8) we have
\[|V_1|^n=|V|\leq 2\sqrt{p}|G|\leq 2\sqrt{p}|G/K||H_1|^n\leq
2\sqrt{p}(\log p)^{8n}|H_1|^n\]
and thus
\[(9)\qquad |H_1|\geq\frac{|V_1|}{2^\frac{1}{n}p^\frac{1}{2n}(\log
p)^8}\]
Next we claim that
\[(10)\qquad |H_1|\geq\frac{|V_1|^\frac{3}{4}}{2}.\]
To see this, first observe that if $n=1$, then $H_1\cong G$ and (10) follows from (6) immediately,
as $m=\dim V\geq 2$.\\
If $n=2$, then we have (by (6) and (8))
\[|V_1|^2\leq 2\sqrt{p}\cdot 2\cdot |H_1|^2\mbox{ and thus }\]
\[|H_1|\geq\frac{|V_1|}{2p^\frac{1}{4}}\geq\frac{|V_1|^\frac{3}{4}}{2}\mbox{,
as claimed.}\]
If $n\geq 3$, then by (9) we have
\[|H_1|\geq\frac{|V_1|}{\sqrt[3]{2}\cdot p^\frac{1}{8}(\log p)^8},\]
which implies (10) if $C$ chosen sufficiently large. This proves
(10).\\
Now we apply \cite[Theorem 3.5(a)]{keller} to the action of $H_1$ on $V_1$ and thus conclude that there are
universal constants $D_1,D_2$ such that with $Z=Z(F^*(H_1))$ we have the following:
\[(11)\qquad\mbox{if }F^*(H_1)=F(H_1)\mbox{, then }|H_1|\leq
D_1|Z|\log |V_1|\mbox{, and}\]
\[(12)\qquad\mbox{if }F^*(H_1)\not=F(H_1)\mbox{, then }|H_1|\leq
D_2|N|\log |V_1|\mbox{, where}\]
$N$ is normal in $H_1$, and
$N=G_1\circ Z$ is a central product of a quasisimple group $G_1$ and $Z$. Moreover, as $Z$ acts fixed point
freely on $|V_1|$, we have
\begin{eqnarray*}
(13)\qquad
k(ZV_1)&=&\frac{|V_1|-1}{|Z|}+|Z|\geq\max\left\{\frac{|V_1|+|Z|-1}{|Z|},|Z|\right\}\\
       &\geq&\sqrt{|V_1|+|Z|-1}\geq|V_1|^\frac{1}{2}
\end{eqnarray*}
Also note that by \ref{cor5} we know that
\[(14)\qquad |H_1|\leq 6|V_1|\log |V_1|.\]
Furthermore, by (7) we have
\[(15)\qquad k(GV)\geq\frac{k(KV)}{|G:K|}\geq\frac{k(KV)}{(\log p)^{8n}},\]
and as $K\ug H_1^n$, we have
\[(16)\qquad
k(KV)\geq\frac{k(H_1^nV)}{|H_1^n:K|}=\frac{k(H_1V_1)^n}{|H_1^n:K|}\]
Now from (6) and (7) we have
\[|K|=\frac{|G|}{|G/K|}\geq\frac{|V_1|^n}{2\sqrt{p}(\log
p)^{8n}},\]
and combining this with (14) yields
\[(17)\qquad |H_1^n:K|\leq 2\sqrt{p}(6(\log p)^8\log
|V_1|)^n\leq\sqrt{p}(7\log |V_1|)^{10n},\]
and this with (15) and (16) shows that
\[(18)\qquad k(GV)\geq\frac{k(H_1V_1)^n}{\sqrt{p}(7\log
|V_1|)^{18n}}=\frac{1}{\sqrt{p}}\left(\frac{k(H_1V_1)}{7^{18}(\log|V_1|)^{18}}\right)^n.\]
$\mbox{ }$\\

Now suppose that $F^*(G)=F(G)$. Then by (11) and (13) we have
\[(19)\qquad
k(H_1V_1)\geq\frac{k(ZV_1)}{|H_1:Z|}\geq\frac{|V_1|^\frac{1}{2}}{D_1\log|V_1|}.\]
Combining (18) and (19) yields
\[(20)\qquad
k(GV)\geq\frac{1}{\sqrt{p}}\left(\frac{|V_1|^\frac{1}{2}}{7^{18}D_1(\log|V_1|)^{19}}\right)^n.\]
If $n=1$, then $H_1=G$, and so if $C$ is large enough, then (19) yields $k(GV)\geq 2\sqrt{p-1}$,
a contradiction, so that we are done in this case. If $n\geq 3$, then by (20) we have
\[(21)\qquad
k(GV)\geq\frac{1}{\sqrt{p}}\frac{|V_1|^\frac{3}{2}}{E(\log|V_1|)^{57}}\geq 2\sqrt{|V_1|-1}\geq 2\sqrt{p-1},\]
for a suitable constant $E$, the first and second inequality being true if $C$ (and thus $|V_1|$) is sufficiently large.
So we have a contradiction as well.\\
It remains to consider the case $n=2$. If $|V_1|\geq p^2$, then by (20) we have, if $C$ is sufficiently large,
\begin{eqnarray*}
k(GV)&\geq&\frac{1}{\sqrt{p}}\frac{|V_1|}{7^{36}D_1^2(\log|V_1|)^{38}}\\
     &\geq&\frac{1}{\sqrt{p}}\frac{p^2}{7^{36}D_1^2(2\log p)^{38}}\\
     &=&\frac{p^\frac{3}{2}}{7^{36}2^{38}D_1^2(\log p)^{38}}\\
     &\geq&2\sqrt{p-1},
\end{eqnarray*}
and we are done.\\
It remains to consider the case $n=2$ and $|V_1|=p$, in which case $H_1$ is cyclic of order dividing $p-1$.
In this case it can be checked by hand that also $k(GV)\geq 2\sqrt{p-1}$, which is the wanted contradiction. Of course,
in this last case, $n=2$ and $|V_1|=p$, we also know $k(GV)\geq 2\sqrt{p-1}$ from \cite{hethelyi}, as $GV$ is
solvable here.\\

So now we may assume that $F^*(H_1)\not=F(H_1)$, thus we have (12). First suppose that $|G_1|\leq C_0$, where
$C_0$ is as in \ref{cor2}.\\
Then similarly as before by (13), (18), (19), and (12) we get
\[(22)\qquad
k(GV)\geq\frac{1}{\sqrt{p}}\left(\frac{|V_1|^\frac{1}{2}}{|H_1:Z|\cdot 7^{18}\cdot(\log|V_1|)^{18}}\right)^n
\geq\frac{1}{\sqrt{p}}\left(\frac{|V_1|^\frac{1}{2}}{7^{18}D_2C_0(\log|V_1|)^{19}}\right)^n\]
which is very similar to (20).\\
If $n=1$, $G\cong H_1$, and so by (19) and (13) we have
\[k(GV)\geq\frac{|V|^\frac{1}{2}}{D_2C_0\log|V|}\]
and thus, as clearly $|V|\geq p^2$, we have
$k(GV)\geq2\sqrt{p-1}$ for $C$ sufficiently large, so that we are done in this case.
So now suppose that $n\geq 2$. Note that as $G_1\leq H_1$, clearly $|V_1|\geq p^2$. As the function
$h(x)=\frac{x^\frac{1}{2}}{D_2C_0(\log x)^{19}}$ is increasing and greater than 1 for large $x$,
we see that if $C$ is large enough,
then (22) implies
\[k(GV)\geq\frac{1}{\sqrt{p}}\left(\frac{\left(p^2\right)^\frac{1}{2}}{7^{18}D_2C_0(\log(p^2))^{19}}\right)^2\geq
\frac{p^\frac{3}{2}}{E_0(\log p)^{19}}\]
for a constant $E_0$, and thus the right hand side is greater than $2\sqrt{p-1}$ if $C$ is large, and we have the
wanted contradiction.\\
So now we may assume that $|G_1|>C_0$. Thus by \ref{cor2} and (12) we have
\[(23)\qquad k(H_1V_1)\geq n(H_1,
V_1)\geq\frac{p}{|H_1:N|}\geq\frac{p}{D_2\log|V_1|},\]
and so by (18) and our assumption $k(GV)<2\sqrt{p-1}$ we obtain
\[2\sqrt{p}>k(GV)\geq\frac{1}{\sqrt{p}}\left(\frac{p}{E_1(\log|V_1|)^{19}}\right)^n\]
and hence
\[(24)\qquad p^{n-1}\leq 2E^n_1(\log|V_1|)^{19n}.\]
Now let $O_\infty(G)$ denote the largest solvable normal subgroup of $G$, and write
$\overline{G}=G/O_\infty(G)$. Clearly $k(\overline{G})\leq k(GV)\leq 2\sqrt{p}$, and by \cite[Lemma 4.7]{pyber}
there is a universal constant $\delta>0$ such that $k(\overline{G})\geq 2^{\delta(\log |\overline{G}|)^\frac{1}{7}}$.
Hence
\[(25)\qquad \log|\overline{G}|\leq\left(\frac{\log
k(\overline{G})}{\delta}\right)^7\leq\left(\frac{1+\frac{1}{2}\log p}{\delta}\right)^7\leq
\left(\frac{\log p}{\delta}\right)^7\]
Moreover, by (8), (7) and (12) we have
\[|G|\leq|G/K||H_1|^n\leq(\log p)^{8n}D_2^n|G_1/G_1\cap
Z|^n|Z|^n(\log|V_1|)^n\]
and as clearly $|G_1/G_1\cap Z|\leq|\overline{G}|$, this together with (25) yields
\[(26)\qquad \log|G|\leq n\left(\log\log p+\log D_2+\left(\frac{\log
p}{\delta}\right)^7+\log|Z|+\log\log|V_1|\right)\]
Therefore, using (6) we see that
\begin{eqnarray*}
(27) n\log|V_1|&=&\log|V|\\
                    &\leq&1+\frac{1}{2}\log p+n\left(\log\log
p+\log D_2+\left(\frac{\log p}{\delta}\right)^7+\log|Z|+\log\log|V_1|\right).
\end{eqnarray*}
Next observe that as $G_1$ is quasisimple, by \cite[Theorem 2.1]{manz/wolf} clearly $(V_1)_Z$
cannot be irreducible and thus is the direct sum of at least two (isomorphic) $Z$--submodules of
$V_1$, and as $Z$ acts fixed point freely on $V_1$, we have $|Z|\leq|V_1|^\frac{1}{2}$.
Therefore $\log |Z|\leq\frac{1}{2}\log|V_1|$. Moreover, if $C$ is sufficiently large, then
$\frac{1}{2}\log|V_1|-\log\log|V_1|\geq\frac{1}{4}\log|V_1|$.
Using all this and (27) implies
\[(28)\qquad\log|V_1|\leq 4\left(1+\log D_2+\frac{1}{2}\log
p+\log\log p+\left(\frac{\log p}{\delta}\right)^7\right)\leq E_2(\log p)^7\]
for a suitable positive constant $E_2$ (as $p>2$).\\
Now if $n=1$, then $G\cong H_1$ and from (23) and (28) we deduce that
\[k(GV)\geq\frac{p}{D_2E_2(\log p)^7}\geq 2\sqrt{p-1},\]
the last inequality following if $C$ is chosen large enough. So this contradiction completes the case
$n=1$.\\
If $n\geq 2$, then (24) implies
\[p\leq 2E_1^\frac{n}{n-1}(\log|V_1|)^\frac{19n}{n-1}\leq
2E_1^2(\log|V_1|)^{38},\]
and thus with (28) we get
\[p\leq 2E_1^2E_2^{19}(\log p)^{266}\]
which for $p>C$ is a contradiction, if $C$ has been chosen sufficiently large. This final contradiction concludes
the proof of the theorem.
\end{beweis}

We can now draw the conclusion that motivated this paper.\\

\begin{ncor}\label{cor8}
Let $C$ be the constant occuring in \ref{theo7}. Let $G$ be a finite group. Suppose $p$ is a prime dividing $|G|$
and that $p>C$. Then
\[k(G)\geq 2\sqrt{p-1}.\]
\end{ncor}
\begin{beweis}
Let $G$ be a counterexample of minimal order. Then by \cite[Section 2]{malle} we know that $G$ has a minimal
normal subgroup $N$ which is an elementary abelian $p$--group, and that $(|G/N|,|N|)=1$. Hence $N$ has a complement
$H\cong G/N$ in $G$, so that $G=HN$, where $N$ is a faithful, irreducible $H$--module over $GF(p)$. Thus by \ref{theo7}
we have $k(G)=k(HN)\geq 2\sqrt{p-1}$ contradicting $G$ being a counterexample. Hence we are done.
\end{beweis}

\end{document}